\documentclass[final,3p, times]{elsarticle}

\usepackage{amssymb,amsmath,color}
\usepackage[mathscr]{eucal}
 \usepackage{lineno}

 \journal{ }
\renewcommand{\journal}{  TBD }
\newtheorem{thm}{Theorem}[section]
\newtheorem{prop}[thm]{Proposition}
\newtheorem{lemma}[thm]{Lemma}
\newtheorem{cor}[thm]{Corollary}

\newproof{pf}{Proof}

\begin{document}

\begin{frontmatter}



\title{$t$-singular linear spaces}

\author[JG]{Jun Guo}\ead{guojun$_-$lf@163.com}
\author[KW]{Kaishun Wang\corref{cor}}
\ead{wangks@bnu.edu.cn}
\author[FL]{Fenggao Li}\ead{fenggaol@163.com}
\cortext[cor]{Corresponding author}

\address[JG]{Math. and Inf. College, Langfang Teachers' College, Langfang  065000,  China }
\address[KW]{Sch. Math. Sci. \& Lab. Math. Com. Sys.,
Beijing Normal University, Beijing  100875, China}
\address[FL]{Dept. of Math., Hunan Institute of Science and Technology, Yueyang  414006,   China}

\begin{abstract}
As a generalization of singular linear spaces, we introduce the concept of
$t$-singular linear spaces, make some anzahl formulas of subspaces,
and determine the suborbits of $t$-singular linear groups.
\end{abstract}

\end{frontmatter}

\section{Introduction}

Let $\mathbb{F}_q$ be a finite field with $q$ elements,
where $q$ is a prime power. For a non-negative integer $n$,
$\mathbb{F}_q^{n}$ denotes the $n$-dimensional row vector space
over $\mathbb{F}_q$.
The set of all $(n_1+n_2)\times (n_1+n_2)$ nonsingular matrices
over $\mathbb{F}_q$
$$\bordermatrix{ &_{n_1}&_{n_2}\cr
&T_{11}&T_{12}\cr
&0&T_{22}}\hspace{-3pt}
\begin{array}{l}
_{n_1}\\
 _{n_2}\end{array}$$
forms a group under matrix multiplication,
called the {\it singular general linear group} of degree $n_1 + n_2$ over $\mathbb{F}_q$
and denoted by $GL_{n_1,n_2}(\mathbb{F}_q)$.
The vector space $\mathbb{F}_{q}^{n_1+n_2}$ together with the right
multiplication action of $GL_{n_1,n_2}(\mathbb{F}_q)$ is called the
$(n_1+n_2)$-dimensional {\it singular linear space} over $\mathbb{F}_{q}$.

As a generalization of attenuated spaces, the concept of singular
linear spaces was firstly introduced in \cite{WGL}.
In \cite{WGL2,WGL3}, we obtained some anzahl formulas of subspaces
in singular linear spaces and discussed their applications to
association schemes. In this paper, we generalize the concept of singular linear spaces to
that of $t$-singular linear spaces. In Section 2, we introduce the concept of
$t$-singular linear spaces and give some  anzahl formulas of subspaces.
In Section 3, we determine the suborbits of $t$-singular linear groups.

\section{$t$-singular linear spaces}

This paper involves partitioned matrices whose entries are
themselves submatrices. For typographical convenience, we often
leave blank the zero submatrices. We write $I^{(r)}$ for the
identity matrix of size $r$, and we omit $r$ if it is clear from the
context.

For non-negative integers $n_1,n_2,\ldots,n_t$,
the set of all nonsingular matrices over $\mathbb{F}_q$
$$
\bordermatrix{ &_{n_1}&_{n_2}& &_{n_t}\cr
& T_{11} & T_{12} & \cdots &T_{1t}\cr
& & T_{22} & \cdots &T_{2t}\cr
& & & \ddots & \vdots\cr
&&  & & T_{tt}}\hspace{-3pt}
\begin{array}{l}
_{n_1}\\
 _{n_2}\\
 \\
 _{n_t}\end{array}
$$
forms a group under matrix
multiplication, called the $t$-{\it singular general linear group} of
degree $n_1+n_2+\cdots+n_t$ over $\mathbb{F}_q$ and denoted by
$GL_{n_1,n_2,\ldots,n_t}(\mathbb{F}_q)$.

Let $P$ be an $m$-dimensional subspace of $\mathbb{F}_q^{n_1+n_2+\cdots+n_t}$,
denote also by $P$ an $m\times (n_1+n_2+\cdots+n_t)$ matrix of rank $m$ whose rows
span the subspace $P$ and call the matrix $P$ a matrix
representation of the subspace $P$. There is an action of
$GL_{n_1,n_2,\ldots,n_t}(\mathbb{F}_q)$ on $\mathbb{F}_q^{n_1+n_2+\cdots+n_t}$ defined as
follows:
\begin{eqnarray}
\mathbb{F}_q^{n_1+n_2+\cdots+n_t}\times GL_{n_1,n_2,\ldots,n_t}(\mathbb{F}_q)&
\longrightarrow&
\mathbb{F}_q^{n_1+n_2+\cdots+n_t}\nonumber\\
((x_1,\ldots,x_{n_1+n_2+\cdots+n_t}),T)&\longmapsto&
(x_1,\ldots,x_{n_1+n_2+\cdots+n_t})T.\nonumber
\end{eqnarray}
The above action induces an action on the set of subspaces of
$\mathbb{F}_q^{n_1+n_2+\cdots+n_t}$; i.e., a subspace $P$ is carried by $T\in
GL_{n_1,n_2,\ldots,n_t}(\mathbb{F}_q)$ to the subspace $PT$. The vector space
$\mathbb{F}_q^{n_1+n_2+\cdots+n_t}$ together with the above group action is called
the $(n_1+n_2+\cdots+n_t)$-dimensional $t$-{\it singular linear space} over
$\mathbb{F}_q$. Note that 2-singular linear spaces are just the singular linear spaces.

For $1\leq j\leq n_1+n_2+\cdots+n_t$, let $e_j $ be the row vector in
$\mathbb{F}_q^{n_1+n_2+\cdots+n_t}$ whose $j$th coordinate is 1 and all other
coordinates are 0. For $2\leq i\leq t$, denote by $E_i$ the $(n_i+n_{i+1}+\cdots+n_t)$-dimensional subspace of
$\mathbb{F}_q^{n_1+n_2+\cdots+n_t}$ generated by $e_{n_1+\cdots+n_{i-1}+1},e_{n_1+\cdots+n_{i-1}+2},\ldots,e_{n_1+n_2+\cdots+n_t}$.
A $k_1$-dimensional subspace $P$ of $\mathbb{F}_q^{n_1+n_2+\cdots+n_t}$ is called a
 subspace of {\it type} $(k_1,k_2,\ldots,k_{t})$ if $\dim(P\cap E_i)=k_i$ for each $2\leq i\leq t$.

Let $\mathcal{M}(k_1,k_2,\ldots,k_{t};n_1,n_2,\ldots,n_t)$ denote the set of all the subspaces of
type $(k_1,k_2,\ldots,k_{t})$ of $\mathbb{F}_q^{n_1+n_2+\cdots+n_t}$  and let
$N(k_1,k_2,\ldots,k_{t};n_1,n_2,\ldots,n_t)$ be the size of $\mathcal{M}(k_1,k_2,\ldots,k_{t};n_1,n_2,\ldots,n_t)$.

\begin{thm}\label{lem2.1}
The set $\mathcal{M}(k_1,k_2,\ldots,k_{t};n_1,n_2,\ldots,n_t)$ is non-empty if and only if
\begin{equation}\label{shi2.1.nzj}
0\leq k_i-k_{i+1}\leq n_i\,(1\leq i\leq t-1)\; \hbox{and}\; 0\leq k_t\leq n_t.
\end{equation} Moreover,
if (\ref{shi2.1.nzj}) holds, then $\mathcal{M}(k_1,k_2,\ldots,k_{t};n_1,n_2,\ldots,n_t)$
forms an orbit under
$GL_{n_1,n_2,\ldots,n_t}(\mathbb{F}_q)$ and
$$
N(k_1,k_2,\ldots,k_{t};n_1,n_2,\ldots,n_t)=\left[n_t\atop k_t\right]_q
\prod\limits_{j=1}^{t-1}q^{(k_j-k_{j+1})(n_{j+1}+\cdots+n_{t}-k_{j+1})}\left[n_j\atop k_j-k_{j+1}\right]_q.
$$
\end{thm}

\begin{pf} The first statement is trivial.

Now suppose (\ref{shi2.1.nzj}) holds.
Pick any subspace $P$ of type $(k_1,k_2,\ldots,k_{t})$
in $\mathbb{F}_q^{n_1+n_2+\cdots+n_t}$. Then $P$ has a matrix representation
\begin{equation}\label{shi2.2nzj}
\bordermatrix{ &_{n_1}&_{n_2}&& _{n_{t-1}}&_{n_t}\cr
&P_{11}&P_{12}&\cdots & P_{1,t-1}& P_{1t}\cr
&&P_{22}& \cdots & P_{2,t-1}& P_{2t}\cr
& &  & \ddots & \vdots & \vdots\cr
& &  & & P_{t-1,t-1}& P_{t-1,t}\cr
&  &  &  & & P_{tt}}\hspace{-3pt}
\begin{array}{l}
_{k_1-k_2}\\
 _{k_2-k_3}\\
\\
 _{k_{t-1}-k_t}\\
 _{k_t}\end{array},
\end{equation}
where  rank $P_{ii}=k_i-k_{i+1}\,(1\leq i\leq t-1)$
and rank $P_{tt}=k_t.$
By basic facts in linear algebra, there exists a
$T\in GL_{n_1,n_2,\ldots,n_t}(\mathbb{F}_q)$ such that
$PT$ has the matrix representation
\begin{equation}\label{shi2.3nzj}
\bordermatrix{ &_{n_1}&_{n_2}&& _{n_{t-1}}&_{n_t}\cr
&(I\;0)&& & & \cr
&&(I\;0)&  & & \cr
&  &  &  \ddots& & \cr
&  &  && (I\;0)& \cr
&  &  & &  & (I\;0)}\hspace{-3pt}
\begin{array}{l}
_{k_1-k_2}\\
 _{k_2-k_3}\\
\\
 _{k_{t-1}-k_t}\\
 _{k_t}\end{array}.
\end{equation}
Hence, $\mathcal{M}(k_1,k_2,\ldots,k_{t};n_1,n_2,\ldots,n_t)$  forms an orbit under
$GL_{n_1,n_2,\ldots,n_t}(\mathbb{F}_q)$.

Denote by $n(k_1,k_2,\ldots,k_{t};n_1,n_2,\ldots,n_t)$ the number of
all matrices of the form (\ref{shi2.2nzj}). Suppose that $P$ and $Q$
represent the same subspace, then there is a  $k_1\times k_1$
nonsingular matrix $U$ such that $P=UQ$. It follows that $U$ is
necessarily of the form
$$
U=\left(\begin{array}{cccc}
 U_{11} & U_{12}&\cdots&U_{1t}\\
        & U_{22}&\cdots&U_{2t}\\
          &&\ddots&\vdots\\
          &&&U_{tt}
\end{array}\right)\in GL_{k_1-k_2,k_2-k_3,\ldots,k_{t-1}-k_t,k_t}(\mathbb{F}_q).
$$
Moreover, if $UQ=Q$, then $U=I$. Consequently,
$$
n(k_1,k_2,\ldots,k_{t};n_1,n_2,\ldots,n_t)=|GL_{k_1-k_2,k_2-k_3,\ldots,k_{t-1}-k_t,k_t}(\mathbb{F}_q)|
N(k_1,k_2,\ldots,k_{t};n_1,n_2,\ldots,n_t).
$$
By \cite[Lemma 1.5]{Wan}, the desired result follows. \qed  \end{pf}

For a fixed subspace $P$ of type $(k_1,k_2,\ldots,k_{t})$
in $\mathbb{F}_q^{n_1+n_2+\cdots+n_t}$,
let $\mathcal{M}(l_1,l_2,\ldots,l_{t};k_1,k_2,\ldots,k_{t};n_1,n_2,\ldots,n_t)$ be the set of all the
subspaces of type $(l_1,l_2,\ldots,l_{t})$ contained in $P$, and let
 $$
 N(l_1,l_2,\ldots,l_{t};k_1,k_2,\ldots,k_{t};n_1,n_2,\ldots,n_t)
 =|\mathcal{M}(l_1,l_2,\ldots,l_{t};k_1,k_2,\ldots,k_{t};n_1,n_2,\ldots,n_t)|.
 $$
  By the transitivity of
$GL_{n_1,n_2,\ldots,n_t}(\mathbb{F}_q)$ on the set of subspaces of the same
type, $N(l_1,l_2,\ldots,l_{t};k_1,k_2,\ldots,k_{t};n_1,n_2,\ldots,n_t)$
is independent of the particular choice
of the subspace $P$ of type $(k_1,k_2,\ldots,k_{t})$.

\begin{prop}\label{lem2.2}
The set $\mathcal{M}(l_1,l_2,\ldots,l_{t};k_1,k_2,\ldots,k_{t};n_1,n_2,\ldots,n_t)$
is non-empty if and only if
\begin{equation}\label{shi1}
0\leq l_i-l_{i+1}\leq k_i-k_{i+1}\leq n_i\,(1\leq i\leq t-1) \;
\hbox{and}\;0\leq l_t\leq k_t\leq n_t.
\end{equation}
Moreover, if (\ref{shi1}) holds, then
\begin{equation}\label{shi2}
N(l_1,l_2,\ldots,l_{t};k_1,k_2,\ldots,k_{t};n_1,n_2,\ldots,n_t)
=\left[k_t\atop l_t\right]_q
\prod\limits_{j=1}^{t-1}q^{(l_j-l_{j+1})(k_{j+1}-l_{j+1})}
\left[k_j-k_{j+1}\atop l_j-l_{j+1}\right]_q.
\end{equation}
\end{prop}

\begin{pf} The first statement is trivial.

Suppose (\ref{shi1}) holds.
By the transitivity of $GL_{n_1,n_2,\ldots,n_t}(\mathbb{F}_q)$ on the
set of subspaces of the same type, we may pick the subspace $P$ of
 type $(k_1,k_2,\ldots,k_{t})$  as the form (\ref{shi2.3nzj}).
Since the number of subspaces of type $(l_1,l_2,\ldots,l_{t})$
in $\mathbb{F}_q^{n_1+n_2+\cdots+n_t}$  contained in $P$ is equal to the number of
subspaces of type $(l_1,l_2,\ldots,l_{t})$ in $\mathbb{F}_q^{(k_1-k_2)+(k_2-k_3)+\cdots+(k_{t-1}-k_t)+k_t}$,
by Lemma~\ref{lem2.1}, (\ref{shi2}) holds. \qed \end{pf}

For a fixed subspace $P$ of type $(l_1,l_2,\ldots,l_{t})$
in $\mathbb{F}_q^{n_1+n_2+\cdots+n_t}$, let $\mathcal{M}'(l_1,l_2,\ldots,l_{t};k_1,k_2,\ldots,k_{t};n_1,n_2,\ldots,n_t)$  be
the set of all the subspaces of type $(k_1,k_2,\ldots,k_{t})$ containing $P$, and let
$$
 N'(l_1,l_2,\ldots,l_{t};k_1,k_2,\ldots,k_{t};n_1,n_2,\ldots,n_t)
 =|\mathcal{M}'(l_1,l_2,\ldots,l_{t};k_1,k_2,\ldots,k_{t};n_1,n_2,\ldots,n_t)|.
 $$
By the transitivity of
$GL_{n_1,n_2,\ldots,n_t}(\mathbb{F}_q)$ on the set of subspaces of the same
type, $N'(l_1,l_2,\ldots,l_{t};k_1,k_2,\ldots,k_{t};n_1,n_2,\ldots,n_t)$ is independent of the particular
choice of the subspace $P$ of type $(l_1,l_2,\ldots,l_{t})$.
By Proposition~\ref{lem2.2}, $\mathcal{M}'(l_1,l_2,\ldots,l_{t};k_1,k_2,$ $\ldots,k_{t};n_1,n_2,\ldots,n_t)$
is non-empty if and only if (\ref{shi1}) holds.

\begin{cor}\label{lem2.3}
If (\ref{shi1}) holds, then
\begin{equation}\label{shi3}
N'(l_1,l_2,\ldots,l_{t};k_1,k_2,\ldots,k_{t};n_1,n_2,\ldots,n_t)
= q^{\sum\limits_{j=1}^{t-1}(k_j-k_{j+1}-l_j+l_{j+1})(n_{j+1}+\cdots+n_t-k_{j+1})}
\left[n_t-l_t\atop k_t-l_t\right]_q
\prod\limits_{j=1}^{t-1}\left[n_j-l_j+l_{j+1}\atop k_j-k_{j+1}-l_j+l_{j+1}\right]_q.
\end{equation}
\end{cor}

\begin{pf} Let
$$
M=\{(P, Q)\mid P\in \mathcal{M}(l_1,l_2,\ldots,l_{t};n_1,n_2,\ldots,n_t), Q\in
\mathcal{M}(k_1,k_2,\ldots,k_{t};n_1,n_2,\ldots,n_t), P\subseteq Q\}.
$$
By computing the size of $M$  in two ways, we have
\begin{eqnarray*}
&&N'(l_1,l_2,\ldots,l_{t};k_1,k_2,\ldots,k_{t};n_1,n_2,\ldots,n_t)
N(l_1,l_2,\ldots,l_{t};n_1,n_2,\ldots,n_t)\\
&=&N(l_1,l_2,\ldots,l_{t};k_1,k_2,\ldots,k_{t};n_1,n_2,\ldots,n_t)
 N(k_1,k_2,\ldots,k_{t};n_1,n_2,\ldots,n_t).
\end{eqnarray*}
By Theorem~\ref{lem2.2}, (\ref{shi3})
holds. \qed \end{pf}

\section{Suborbits}

Let $G$ be a transitive permutation group on a finite set $\Omega$,
denoted by $(G,\Omega)$.
For a fixed element $a\in\Omega$, the orbits of $G_a$
on $\Omega$ are said to be the {\it suborbits} of $(G,\Omega)$, and the number of such
suborbits is the {\it rank} of $(G,\Omega)$. The size of each suborbit
is said to be the its {\it length}.
The results on suborbits of classical groups on the set of subspaces may be found in Wang and Wei (\cite{Wang}),
Wei and Wang (\cite{Wei1,Wei2}), Guo, Wang and Li \cite{GUO,gw,GWL2,GWL3}).

Let $$U=\bordermatrix{ &_{n_1}&_{n_2}&_{n_3}\cr
&U_{11}&&\cr
&&U_{22}&\cr
&&&U_{33}}\hspace{-3pt}
\begin{array}{c}
_{k_1-k_2}\\
_{k_2-k_3}\\
_{k_3}
\end{array}=\bordermatrix{ &_{n_1}&_{n_2}&_{n_3}\cr
&(I\;0)&&\cr
&&(I\;0)&\cr
&&&(I\;0)}\hspace{-3pt}
\begin{array}{c}
_{k_1-k_2}\\
_{k_2-k_3}\\
_{k_3}
\end{array}\in\mathcal{M}(k_1,k_2,k_3;n_1,n_2,n_3)$$
and let $G_U$ be the stabilizer of $U$ in $GL_{n_1,n_2,n_3}(\mathbb{F}_q)$.
In order to determine the suborbits of
$(GL_{n_1,n_2,n_3}(\mathbb{F}_q),\mathcal{M}(k_1,k_2,k_3;n_1,n_2,n_3))$,
we only need to consider the orbits of $G_U$ on $\mathcal{M}(k_1,k_2,k_3;n_1,n_2,n_3)$.

\begin{thm}\label{thm3.1}
Let $0\leq k_3\leq n_3,0\leq k_2-k_3\leq n_2$ and $0\leq k_1-k_2\leq n_1$.
Two elements of
$\mathcal{M}(k_1,k_2,k_3;n_1,n_2,n_3)$
$$
Q=\bordermatrix{ &_{n_1}&_{n_2}&_{n_3}\cr
&Q_{11}&Q_{12}&Q_{13}\cr
&&Q_{22}&Q_{23}\cr
&&&Q_{33}}\hspace{-3pt}
\begin{array}{c}
_{k_1-k_2}\\
_{k_2-k_3}\\
_{k_3}
\end{array},\;P=\bordermatrix{ &_{n_1}&_{n_2}&_{n_3}\cr
&P_{11}&P_{12}&P_{13}\cr
&&P_{22}&P_{23}\cr
&&&P_{33}}\hspace{-3pt}
\begin{array}{c}
_{k_1-k_2}\\
_{k_2-k_3}\\
_{k_3}
\end{array},$$
fall  into the same orbit of
$G_U$ if and only if
\begin{equation}\label{s0}
\left.\begin{array}{l}\dim(U_{ii}\cap Q_{ii})=\dim(U_{ii}\cap P_{ii}) \;(i=1,2,3),\;
\dim(U\cap Q)=\dim(U\cap P),\\
\dim\left(\left(\begin{array} {cc}U_{11}&0\\
0& U_{22} \end{array} \right)
\cap\left(\begin{array} {cc}Q_{11}&Q_{12}\\
0& Q_{22} \end{array} \right)\right)
=\dim\left(\left(\begin{array} {cc}
U_{11}&0\\
0&U_{22} \end{array} \right)
\cap\left(\begin{array} {cc}
P_{11}&P_{12}\\
0&P_{22} \end{array} \right)\right),\\
\dim\left(\left(\begin{array} {cc}U_{22}&0\\
0& U_{33} \end{array} \right)
\cap \left(\begin{array} {cc}Q_{22}&Q_{23}\\
0& Q_{33} \end{array} \right)\right)
=\dim\left(\left(\begin{array} {cc}
U_{22}&0\\
0&U_{33} \end{array} \right)
\cap\left(\begin{array} {cc}
P_{22}&P_{23}\\
0&P_{33} \end{array} \right)\right).\end{array}\right\}\end{equation}
\end{thm}

\begin{pf} Suppose $Q$ and $P$ are in the same orbit of
$G_U$. Then there exists
a
$$T=\left(
\begin{array} {ccc}
  T_{11}& T_{12}&T_{13}\\
 & T_{22}&T_{23}\\
 &&T_{33}
 \end{array} \right)\in G_U$$
  such that
\begin{eqnarray*}
QT&=&\left(
\begin{array} {ccc}
Q_{11}T_{11}&
Q_{11}T_{12}+Q_{12}T_{22}&Q_{11}T_{13}+Q_{12}T_{23}+Q_{13}T_{33}\\
& Q_{22}T_{22}&Q_{22}T_{23}+Q_{23}T_{33}\\
&&Q_{33}T_{33}\end{array} \right)=\left(
\begin{array} {ccc}
P_{11}&P_{12}&P_{13}\\
&P_{22}&P_{23}\\
&&P_{33} \end{array} \right)=P.
\end{eqnarray*}
Hence (\ref{s0}) holds.

Conversely, suppose (\ref{s0}) holds. Let
\begin{equation}\label{s01}
\left.\begin{array}{l}
  \dim(U_{11}\cap Q_{11})=k_1-k_2-i_1,
 \dim(U_{22}\cap Q_{22})=k_2-k_3-i_2,\\
 \dim(U_{33}\cap Q_{33})=k_3-i_3,
 \dim(U\cap Q)=k_1-j_1,\\
 \dim\left(\left(\begin{array} {cc}U_{11}&0\\
0& U_{22} \end{array} \right)
\cap\left(\begin{array} {cc}Q_{11}&Q_{12}\\
0& Q_{22} \end{array} \right)\right)
=k_1-k_3-j_{2},\\
\dim\left(\left(\begin{array} {cc}U_{22}&0\\
0& U_{33} \end{array} \right)
\cap \left(\begin{array} {cc}Q_{22}&Q_{23}\\
0& Q_{33} \end{array} \right)\right)
=k_2-j_{3}.\end{array}\right\}.
\end{equation}
Then $U$ and $Q$ have the matrix representations
\begin{equation}\label{s1} U=\bordermatrix{&_{n_1}&_{n_2}&_{n_3}\cr
 & U_{111}&&\cr
  &U_{112}&&\cr
  &U_{113}&&\cr
   &U_{114}&&\cr
  &&U_{221}&\cr
  &&U_{222}&\cr
  &&U_{223}&\cr
  &&&U_{331}\cr
  &&&U_{332}}\hspace{-3pt}
\begin{array}{c}
_{i_1}\\
_{i_2+j_1-j_{2}-j_{3}}\\
_{k_1-k_2+j_{3}-j_1}\\
_{j_{2}-i_1-i_2}\\
_{i_2}\\
_{k_2-k_3-j_{3}+i_3}\\
_{j_{3}-i_2-i_3}\\
_{i_3}\\
_{k_3-i_3}
\end{array} \quad\mbox{and}\quad
Q=\bordermatrix{&_{n_1}&_{n_2}&_{n_3}\cr
 & Q_{111}&Q_{121}&Q_{131}\cr
  &U_{112}&0&Q_{132}\cr
  &U_{113}&0&0\cr
   &U_{114}&Q_{124}&Q_{134}\cr
  &&Q_{221}&Q_{231}\cr
  &&U_{222}&0\cr
  &&U_{223}&Q_{233}\cr
  &&&Q_{331}\cr
  &&&U_{332}}\hspace{-3pt}
\begin{array}{c}
_{i_1}\\
_{i_2+j_1-j_{2}-j_{3}}\\
_{k_1-k_2+j_{3}-j_1}\\
_{j_{2}-i_1-i_2}\\
_{i_2}\\
_{k_2-k_3-j_{3}+i_3}\\
_{j_{3}-i_2-i_3}\\
_{i_3}\\
_{k_3-i_3}
\end{array},
\end{equation}
where rank $Q_{233}=j_{3}-i_2-i_3$, rank $Q_{124}=j_{2}-i_1-i_2$
and rank $Q_{132}=i_2+j_1-j_{2}-j_{3}$.
It follows that $U+Q$ is a
subspace of type $(k_1+j_1,k_2+j_{1}-i_1,k_3+j_1-j_2)$ with a matrix representation of the
form
$$
\left(\begin{array}{ccc}
  U_{111}&&\\
  U_{112}&&\\
  U_{113}&&\\
   U_{114}&&\\
   Q_{111}&Q_{121}&Q_{131}\\
   &-Q_{124}&-Q_{134}\\
  &U_{221}&0\\
  &U_{222}&0\\
  &U_{223}&0\\
  &Q_{221}&Q_{231}\\
  &&-Q_{132}\\
  &&-Q_{233}\\
  &&U_{331}\\
  &&U_{332}\\
  &&Q_{331}
\end{array}\right).
$$
Similarly,  $U+P$ is also a subspace of type $(k_1+j_1,k_2+j_{1}-i_1,k_3+j_1-j_2)$ with a
matrix representation just like that of $U+Q$. By Theorem~\ref{lem2.1}, there exists a
$T\in GL_{n_1,n_2,n_3}(\mathbb{F}_q)$ such that $(U+P)T=U+Q$, which
implies that $UT=U$ and $PT=Q$. Hence both $Q$ and $P$ are
in the same orbit of $G_U$. \qed \end{pf}

For any $Q$ of the form (\ref{s1}),
let $\Lambda_{(i_1,i_2,i_3,j_3-i_2-i_3,j_2-i_1-i_2,i_2+j_1-j_2-j_3)}$
be the orbit of $G_U$ containing $Q$. Then
\begin{eqnarray}
&&0\leq i_1\leq \min\{k_1-k_2,n_1+k_2-k_1\},\;
0\leq i_2\leq \min\{k_2-k_3,n_2+k_3-k_2\},\;
0\leq i_3\leq \min\{k_3,n_3-k_3\},\label{shin9}\\
&&\max\{k_2-k_3-i_2,k_3-i_3\}\leq k_2-j_3\leq (k_2-k_3-i_2)+(k_3-i_3),\label{shin10}\\
&&\max\{k_1-k_2-i_1,k_2-k_3-i_2\}\leq k_1-k_3-j_2\leq (k_1-k_2-i_1)+(k_2-k_3-i_2),\label{shin11}\\
&&k_2-j_3\leq k_1-j_1\leq (k_1-k_3-j_2)+(k_2-j_3)-(k_2-k_3-i_2),\label{shin12}\\
&&k_3+j_1-j_2\leq n_3,\,
k_2-k_3+j_2-i_1\leq n_2.\label{shin13}
\end{eqnarray}
By (\ref{shin10})-(\ref{shin13}),
\begin{eqnarray*}
&&0\leq j_3-i_2-i_3\leq \min\{k_3-i_3,k_2-k_3-i_2\},\\
&&0\leq j_2-i_1-i_2\leq \min\{k_2-k_3-i_2,k_1-k_2-i_1,n_2+k_3-k_2-i_2\},\\
&&0\leq i_2+j_1-j_2-j_3\leq \min\{k_1-k_2+i_2-j_2,n_3-k_3+i_2-j_3\}.
\end{eqnarray*}
Then we have
\begin{equation}\label{shin14}
\left.\begin{array}{l}
0\leq i_1\leq \min\{k_1-k_2,n_1+k_2-k_1\},\\
0\leq i_2\leq \min\{k_2-k_3,n_2+k_3-k_2\},\\
0\leq i_3\leq \min\{k_3,n_3-k_3\},\\
0\leq j_3-i_2-i_3\leq \min\{k_3-i_3,k_2-k_3-i_2\},\\
0\leq j_2-i_1-i_2\leq \min\{k_2-k_3-i_2,k_1-k_2-i_1,n_2+k_3-k_2-i_2\},\\
0\leq i_2+j_1-j_2-j_3\leq \min\{k_1-k_2+i_2-j_2,n_3-k_3+i_2-j_3\}.
\end{array}\right\}
\end{equation}

Conversely, for any given integers $i_1,i_2,i_3,j_1,j_2$ and $j_3$ satisfying
(\ref{shin14}), pick
$$
Q=\bordermatrix{&_{i_1}&_{k_1-k_2-i_1}&_{i_1}&_{n_1-k_1+k_2-i_1}&_{i_2}&_{k_2-k_3-i_2}
&_{i_2}&_{n_2-k_2+k_3-i_2}&_{i_3}&_{k_3-i_3}&_{i_3}&_{n_3-k_3-i_3}\cr
 & 0&I&0&0&0&0&0&A&0&0&0&C\cr
 & 0&0&I&0&0&0&0&0&0&0&0&0\cr
 & 0&0&0&0&0&I&0&0&0&0&0&B\cr
 & 0&0&0&0&0&0&I&0&0&0&0&0\cr
 & 0&0&0&0&0&0&0&0&0&I&0&0\cr
 & 0&0&0&0&0&0&0&0&0&0&I&0}
\begin{array}{c}
_{k_1-k_2-i_1}\\
_{i_1}\\
_{k_2-k_3-i_2}\\
_{i_2}\\
_{k_3-i_3}\\
 _{i_3}
\end{array},
$$
where
$$A=\left(\begin{array}{cc}
0&0\\
I^{(j_2-i_1-i_2)}&0
\end{array}\right),\;
B=\left(\begin{array}{cc}
0&0\\
I^{(j_3-i_2-i_3)}&0
\end{array}\right),\;
C=\left(\begin{array}{cc}
0&I^{(i_2+j_1-j_2-j_3)}\\
0&0
\end{array}\right).$$
Then $Q\in \Lambda_{(i_1,i_2,i_3,j_3-i_2-i_3,j_2-i_1-i_2,i_2+j_1-j_2-j_3)}$; and so the orbit
$\Lambda_{(i_1,i_2,i_3,j_3-i_2-i_3,j_2-i_1-i_2,i_2+j_1-j_2-j_3)}$ exists. It follows that the orbits of
$G_U$ are completely
determined by $(i_1,i_2,i_3,j_3-i_2-i_3,j_2-i_1-i_2,i_2+j_1-j_2-j_3)$ satisfying (\ref{shin14}).
Therefore, we have the following result.

\begin{thm}\label{thm2} Let $0\leq k_1-k_2\leq n_1,0\leq k_2-k_3\leq n_2$ and $0\leq k_3\leq n_3$.
Then the number of  suborbits of
$(GL_{n_1,n_2,n_3}(\mathbb{F}_q),\mathcal{M}(k_1,k_2,k_3;n_1,n_2,n_3))$ is
$$\sum\limits_{i_1=0}^{\min\{k_1-k_2,n_1+k_2-k_1\}}\sum\limits_{i_2=0}^{\min\{k_2-k_3,n_2+k_3-k_2\}}
\sum\limits_{i_3=0}^{\min\{k_3,n_3-k_3\}}\sum\limits_{j_3=i_2+i_3}^{\min\{k_3+i_2,k_2-k_3+i_3\}}
\sum\limits_{j_2=i_1+i_2}^{\min\{k_2-k_3+i_1,k_1-k_2+i_2,n_2+k_3-k_2+i_1\}}
(1+\min\{k_1-k_2+i_2-j_2,n_3-k_3+i_2-j_3\}).$$
\end{thm}

\medskip
In order to compute the length of suborbits of
$(GL_{n_1,n_2,n_3}(\mathbb{F}_q),\mathcal{M}(k_1,k_2,k_3;n_1,n_2,n_3))$, we need the
following   results.

\begin{prop}\label{prop2.1} {\rm(\cite[Chapter 1, Theorem~5]{Wan2})}
The number of $m\times n$ matrices with rank $i$ over $\mathbb{F}_q$
is
$$
N(i;m\times n)=q^{i(i-1)/2}\,\left[m\atop
i\right]_q\prod\limits^{n}_{t=n-i+1}(q^t-1).
$$
\end{prop}

\begin{prop}\label{prop2.2}{\rm(\cite[Chapter 6, Theorem~7]{Wan2})}
Let $1\leq m\leq n$ and $0\leq i\leq\min\{m,n-m\}$.
For a given $m$-dimensional subspace $P$ of
$\mathbb{F}_q^{n}$, the number of $m$-dimensional subspaces
intersecting $P$ at $(m-i)$-dimensional subspaces of
$\mathbb{F}_q^{n}$ is
$$
q^{i^2}\left[n-m\atop i\right]_q\left[m\atop i\right]_q.
$$
\end{prop}

\begin{prop}\label{prop2.3} {\rm(\cite[Lemma~2.4]{WGL})}
For any $m_1\times n$ matrix $A_1$ with rank
$t_1$, the number of $m_2\times n$ matrix $A_2$ satisfying
{\rm rank} $\left(A_1\atop A_2\right)=t_2$ is
$q^{m_2t_1} N(t_2-t_1;m_2\times(n-t_1))$.
\end{prop}

\begin{lemma}\label{lem2.4.0} For any $m\times n_1$ matrix $A_1$ with rank
$t_1$, the number of $m\times n_2$ matrix $A_2$ satisfying
{\rm rank} $(A_1\; A_2)=t_2$ is $q^{t_1n_2}N(t_2-t_1;(m-t_1)\times n_2)$.
\end{lemma}

\begin{pf}  The proof is similar to that of \cite[Lemma~2.4]{WGL}, and will be omitted.\qed \end{pf}

\begin{lemma}\label{lem2.4}
The number of $(m_1+m_2)\times (n_1+n_2)$ matrix
$$\bordermatrix{&_{n_1}&_{n_2}\cr
&A&B\cr
&C&D}
\begin{array}{c}
_{m_1}\\
 _{m_2}
\end{array}$$ with {\rm rank} $(C\;D)=\alpha$
and {\rm rank} $\left(B\atop D\right)=\alpha$
 is $$\sum\limits_{l=\max\{0,\alpha-n_1,\alpha-m_1\}}^{\alpha}q^{(m_1+n_1)l+m_1n_1}N(l;m_2\times n_2)N(\alpha-l;m_1\times(n_2-l))
 N(\alpha-l;(m_2-l)\times n_1).$$
\end{lemma}

\begin{pf}
Let rank $D=l$. Then
$\max\{0,\alpha-n_1,\alpha-m_1\}\leq l\leq\alpha$.
By Proposition~\ref{prop2.1}, there are $N(l;m_2\times n_2)$
choices for $D$. For a given $D$,
by Proposition~\ref{prop2.3}
there are $q^{m_1l} N(\alpha-l;m_1\times(n_2-l))$
choices for $B$, by Lemma~\ref{lem2.4.0} there are $q^{ln_1}N(\alpha-l;(m_2-l)\times n_1)$
choices for $C$. Note that there are $q^{m_1n_1}$ choices for $A$. Therefore, the desired
result follows. \qed \end{pf}

\begin{thm}\label{thm3}
Suppose (\ref{shin14}) holds. Then the length of the suborbit
$\Lambda_{(i_1,i_2,i_3,j_3-i_2-i_3,j_2-i_1-i_2,i_2+j_1-j_2-j_3)}$
of $(GL_{n_1,n_2,n_3}(\mathbb{F}_q),$ $\mathcal{M}(k_1,k_2,k_3;n_1,n_2,n_3))$ is
$$q^{(n_2+n_3-k_2)i_1+(n_3+k_1-k_2-k_3-i_1)i_2+(k_1-k_3-i_1-i_2)i_3+i_1^2+i_2^2+i_3^2}\left[n_1+k_2-k_1\atop i_1\right]_q\left[k_1-k_2\atop
i_1\right]_q\left[n_2+k_3-k_2\atop i_2\right]_q\left[k_2-k_3\atop
i_2\right]_q\left[n_3-k_3\atop i_3\right]_q\left[k_3\atop
i_3\right]_q$$
$$\times N(j_3-i_2-i_3;(k_2-k_3-i_2)\times(n_3-k_3-i_3))N(j_2-i_1-i_2;(k_1-k_2-i_1)\times(n_2+k_3-k_2-i_2))$$
$$\times\sum\limits_{l=\max\{0,2i_2+j_1-j_2-2j_3+i_3,2i_2+j_1-2j_2-j_3+i_1\}}^{i_2+j_1-j_2-j_3}
q^{(j_2+j_3-i_1-i_3)l+(j_2-i_1-i_2)(j_3-i_2-i_3)}
N(l;(k_1-k_2-j_2+i_2)\times (n_3-k_3-j_3+i_2))$$
$$\times N(i_2+j_1-j_2-j_3-l;(j_2-i_1-i_2)\times(n_3-k_3-j_3+i_2-l))
N(i_2+j_1-j_2-j_3-l;(k_1-k_2-j_2+i_2-l)\times (j_3-i_2-i_3)).
$$
\end{thm}

\begin{pf} Suppose
$$P=\bordermatrix{ &_{n_1}&_{n_2}&_{n_3}\cr
&P_{11}&P_{12}&P_{13}\cr
&0&P_{22}&P_{23}\cr
&0&0&P_{33}}\hspace{-3pt}
\begin{array}{c}
_{k_1-k_2}\\
_{k_2-k_3}\\
_{k_3}
\end{array}\in\Lambda_{(i_1,i_2,i_3,j_3-i_2-i_3,j_2-i_1-i_2,i_2+j_1-j_2-j_3)}.$$
Then $P_{11}$ is a  $(k_1-k_2)$-dimensional subspace of
$\mathbb{F}^{n_1}_q$ such that $\dim(P_{11}\cap
U_{11})=k_1-k_2-i_1$, $P_{22}$ is a $(k_2-k_3)$-dimensional subspace
of $\mathbb{F}^{n_2}_q$ such that $\dim(P_{22}\cap
U_{22})=k_2-k_3-i_2$ and $P_{33}$ is a $k_3$-dimensional subspace of
$\mathbb{F}^{n_3}_q$ such that $\dim(P_{33}\cap U_{33})=k_3-i_3$. By
Proposition~\ref{prop2.2}, there are
 \begin{eqnarray*}
\alpha=q^{i_1^2+i_2^2+i_3^2}\left[n_1+k_2-k_1\atop i_1\right]_q\left[k_1-k_2\atop
i_1\right]_q\left[n_2+k_3-k_2\atop i_2\right]_q\left[k_2-k_3\atop
i_2\right]_q\left[n_3-k_3\atop i_3\right]_q\left[k_3\atop
i_3\right]_q\end{eqnarray*} choices for $(P_{11},P_{22},P_{33})$. By the
transitivity of $G_U$ on $\Lambda_{(i_1,i_2,i_3,j_3-i_2-i_3,j_2-i_1-i_2,i_2+j_1-j_2-j_3)}$, we may pick
$$P_{11}=(0^{(k_1-k_2,i_1)}\;I^{(k_1-k_2)}\;0^{(k_2-k_1,n_1+k_2-k_1-i_1)}),\;
P_{22}=(0^{(k_2-k_3,i_2)}\;I^{(k_2-k_3)}\;0^{(k_2-k_3,n_2+k_3-k_2-i_2)})
\;\mbox{and}\;P_{33}=(0^{(k_3,i_3)}\;I^{(k_3)}\;0^{(k_3,n_3-k_3-i_3)}).$$
 Then $P_{12},P_{23},P_{13}$ have the matrix representations of the forms
$$
 P_{12}=\bordermatrix{&_{i_2}&_{k_2-k_3}&_{n_2+k_3-k_2-i_2}\cr
 &A_{11}&0&A_{12}\cr
 &A_{21}&0&A_{22}\cr }\hspace{-3pt}
\begin{array}{c}
_{k_1-k_2-i_1}\\
_{i_1}
\end{array},\;
P_{23}=\bordermatrix{&_{i_3}&_{k_3}&_{n_3-k_3-i_3}\cr
 &B_{11}&0&B_{12}\cr
 &B_{21}&0&B_{22}\cr }\hspace{-3pt}
\begin{array}{c}
_{k_2-k_3-i_2}\\
_{i_2}
\end{array},\;
P_{13}=\bordermatrix{&_{i_3}&_{k_3}&_{n_3-k_3-i_3}\cr
 &C_{11}&0&C_{12}\cr
 &C_{21}&0&C_{22}\cr }\hspace{-3pt}
\begin{array}{c}
_{k_1-k_2-i_1}\\
_{i_1}
\end{array},
 $$
where rank $A_{12}=j_2-i_1-i_2$, rank $B_{12}=j_3-i_2-i_3$,
rank $(A_{12}\;C_{12})=j_1-j_3-i_1$ and rank $\left(C_{12}\atop B_{12}\right)=j_1-j_2-i_3$.
By Proposition~\ref{prop2.1}, there are
$N(j_3-i_2-i_3;(k_2-k_3-i_2)\times(n_3-k_3-i_3))$ choices for $B_{12}$,
and there are
$N(j_2-i_1-i_2;(k_1-k_2-i_1)\times(n_2+k_3-k_2-i_2))$ choices for $A_{12}$.

Now we compute the numbers of $C_{12}$ satisfying the above conditions. Let $M$ be the set of all matrices of the forms
\begin{equation}\label{shin15}\bordermatrix{&_{n_2+k_3-k_2-i_2}&_{n_3-k_3-i_3}\cr
 &A_{12}&C_{12}\cr
 &0&B_{12}\cr }\hspace{-3pt}
\begin{array}{c}
_{k_1-k_2-i_1}\\
_{k_2-k_3-i_2}
\end{array},
\end{equation}
where rank $A_{12}=j_2-i_1-i_2$, rank $B_{12}=j_3-i_2-i_3$,
rank $(A_{12}\;C_{12})=j_1-j_3-i_1$ and rank $\left(C_{12}\atop B_{12}\right)=j_1-j_2-i_3$.
Let ${\cal G}$ (resp. ${\cal S}$) be the set of all non-singular matrices of the forms
$$\bordermatrix{&_{k_1-k_2-i_1}&_{k_2-k_3-i_2}\cr
 &T_{11}&\cr
 &&T_{22}\cr }\hspace{-3pt}
\begin{array}{c}
_{k_1-k_2-i_1}\\
_{k_2-k_3-i_2}
\end{array}\;
\left(\hbox{resp.}\;
\bordermatrix{&_{n_2+k_3-k_2-i_2}&_{n_3-k_3-i_3}\cr
 &S_{11}&\cr
 &&S_{22}\cr }\hspace{-3pt}
\begin{array}{c}
_{n_2+k_3-k_2-i_2}\\
_{n_3-k_3-i_3}
\end{array}\right).$$
There is an action of ${\cal G}\times {\cal S}$ on $M$ defined as
follows:
\begin{eqnarray} M\times({\cal G}\times {\cal S}))& \longrightarrow&M\nonumber\\
(A,(T,S))&\longmapsto& T^{-1}AS.\nonumber
\end{eqnarray}
For a given $(A_{12},B_{12})$,   there exist
$$\bordermatrix{&_{k_1-k_2-i_1}&_{k_2-k_3-i_2}\cr
 &T_{11}&\cr
 &&T_{22}\cr }\hspace{-3pt}
\begin{array}{c}
_{k_1-k_2-i_1}\\
_{k_2-k_3-i_2}
\end{array}\in{\cal G}\;
\hbox{and}\;
\bordermatrix{&_{n_2+k_3-k_2-i_2}&_{n_3-k_3-i_3}\cr
 &S_{11}&\cr
 &&S_{22}\cr }\hspace{-3pt}
\begin{array}{c}
_{n_2+k_3-k_2-i_2}\\
_{n_3-k_3-i_3}
\end{array}\in{\cal S}$$
such that
\begin{equation}\label{shin16}\left(\begin{array}{cc}
T_{11}&\\
&T_{22}\end{array}\right)
\left(\begin{array}{cc}
A_{12}&C_{12}\\
&B_{12}\end{array}\right)
\left(\begin{array}{cc}
S_{11}&\\
&S_{22}\end{array}\right)
=\bordermatrix{&_{j_2-i_1-i_2}&_{n_2+k_3-k_2-j_2+i_1}&_{j_3-i_2-i_3}&_{n_3-k_3-j_3+i_2}\cr
 &I&0&C_1&C_2\cr
 &&0&C_3&C_4\cr
 &&&I&0\cr
 &&&&0 }\hspace{-3pt}
\begin{array}{c}
_{j_2-i_1-i_2}\\
_{k_1-k_2-j_2+i_2}\\
_{j_3-i_2-i_3}\\
_{k_2-k_3-j_3+i_3}
\end{array},\end{equation}
where $$\left(\begin{array}{cc}
C_1&C_2\\
C_3&C_4\end{array}\right)=T_{11}C_{12}S_{22},\;{\rm rank}\;(C_3\;C_4)={\rm rank}\;
\left(C_2\atop C_4\right)=i_2+j_1-j_2-j_3.$$
Therefore, for a given $(A_{12},B_{12})$,  the number of
$C_{12}$ satisfying (\ref{shin15}) is equal to the number of
$\left(\begin{array}{cc}
C_1&C_2\\
C_3&C_4\end{array}\right)$ satisfying (\ref{shin16}).
By Lemma~\ref{lem2.4}, the number of
$\left(\begin{array}{cc}
C_1&C_2\\
C_3&C_4\end{array}\right)$ satisfying (\ref{shin16}) is
$$\sum\limits_{l=\max\{0,2i_2+j_1-j_2-2j_3+i_3,2i_2+j_1-2j_2-j_3+i_1\}}^{i_2+j_1-j_2-j_3}
q^{(j_2+j_3-i_1-i_3)l+(j_2-i_1-i_2)(j_3-i_2-i_3)}
N(l;(k_1-k_2-j_2+i_2)\times (n_3-k_3-j_3+i_2))$$
$$\times N(i_2+j_1-j_2-j_3-l;(j_2-i_1-i_2)\times(n_3-k_3-j_3+i_2-l))
N(i_2+j_1-j_2-j_3-l;(k_1-k_2-j_2+i_2-l)\times (j_3-i_2-i_3)).
$$

Note that there are $q^{(n_2+n_3-k_2)i_1+(n_3+k_1-k_2-k_3-i_1)i_2+(k_1-k_3-i_1-i_2)i_3}$
choices for $(A_{11},A_{21},A_{22},B_{11},B_{21},B_{22},C_{11},C_{21},C_{22})$.
Hence the desired result follows. \qed \end{pf}

\noindent{\bf Remarks}.

(i) Dam and Koolen \cite{Dam} constructed the twisted Grassmann graph
$\tilde{J}_q(2e+1,e)$, which is the first know family of non-vertex-transitive distance-regular graphs
with unbounded diameter. Pick the hyperplane $H=(0^{(1,2e)}\;I^{(2e)})$ in $\mathbb{F}_q^{2e+1}$. Then $P_0=(0^{(e-1,e+2)}\;I^{(e-1)})$
is a vertex of the twisted Grassmann graph. Note that
the last subconstituent of $\tilde{J}_q(2e+1,e)$ about $P_0$
is just ${\cal M}(e+1,e,0;1,e+1,e-1)$.

(ii) Similarly, we may determine all the suborbits of
$(GL_{n_1,n_2,\ldots,n_t}(\mathbb{F}_q),\mathcal{M}(k_1,k_2,\ldots,k_t;n_1,n_2,\ldots,n_t))$
for any $t\geq4$.

\section*{Acknowledgement}
This research is partially supported by    NSF of China (10971052,
10871027),   NCET-08-0052, Langfang Teachers' College (LSZB201005),
and  Hunan Provincial Natural Science Foundation
of China (09JJ3006).

\end{document}